\documentclass[notitlepage, 10pt]{article}
\usepackage{amssymb, amsthm, amsmath, amscd, hyperref, graphicx}
\hypersetup{pdfstartview=FitH,pdfauthor={HC}, pdftitle={Newton's Method as a Formal Recurrence}}

\newtheorem{theorem}{Theorem}
\newtheorem{lemma}{Lemma}

\newcommand{\R}{\mathbb{R}}
\newcommand{\C}{\mathbb{C}}
\newcommand{\Cstar}{\mathbb{C}^*}
\newcommand{\onto}{\rightarrow}
\newcommand{\qbinomial}[2]{\genfrac{[}{]}{0pt}{}{ #1 }{ #2 }_q }

\title{Newton's Method as a Formal Recurrence}
\author{Hal Canary, Carl Edquist, Samuel Lachterman, Brendan Younger}

\begin{document}

\maketitle

\begin{abstract}
Iterating Newton's method symbolically for the general quadratic $ ax^2+bx+c $ yields a rational function $ \frac{P_n(x)}{Q_n(x)} $, the numerator and denominator of which are polynomials with highly composite coefficients. In particular, the coefficients have no prime factors greater than $ 2^n $ after $ n $ iterations even though most of the coefficients are much larger than $ 2^n $.
\end{abstract}


\section{Introduction}

If $f:\C \onto \C$ is a differentiable function, then Newton's method, applied recursively to an initial value of $ z_0 $, yields the sequence of values $ z_1, z_2, \ldots  $ defined as
\begin{equation}
z_{n+1} = N(z_n) = z_n - \frac{f(z_n)}{f'(z_n)}
\label{eq:newtonrecurrence}
\end{equation}
which, in many cases, converges to a root of $ f $. We restrict our attention to the general quadratic $ f(x) = ax^2 + bx + c $ for the purposes of this paper. Instead of using equation (\ref{eq:newtonrecurrence}) as a numerical method, we are interested in iterating symbolically by letting the initial guess be $ z_0 = x $ where $ x \in \C $.  Doing so, we obtain
$$ z_n = \frac{P_n(x)}{Q_n(x)} $$
where $ P_n(x) $ and $ Q_n(x) $ are relatively prime polynomials in $ \C[x] $. An interesting observation, pointed out to us by Jim Propp, is that the coefficients of these polynomials, though very large, have no prime factors $ \geq 2^n $.  On his suggestion that this divisibility property might imply some combinatorial interpretation not immediately apparent in the formulation of equation (\ref{eq:newtonrecurrence}), we have derived an explicit symbolic formula for each iterate.

The aim of this paper is to show that the polynomials $ P_n $ and $ Q_n $ are given by the following explicit formulas:
\begin{subequations}
\label{eq:explicitpolynomials}
\begin{equation}
P_n(x) = a^{2^n-1} x^{2^n} - \sum_{k=0}^{2^n-2} ~ \sum_{j=0}^{2^n-k-2} (-1)^j \binom{2^n}{k} \binom{2^n-k-j-2}{j} a^{k+j} b^{2^n-k-2j-2} c^{j+1} x^k 
\label{eq:explicitp}
\end{equation}
\begin{equation}
Q_n(x) = \sum_{k=0}^{2^n-1}~ \sum_{j=0}^{2^n-k-1} (-1)^j \binom{2^n}{k} \binom{2^n-k-j-1}{j} a^{k+j} b^{2^n-k-2j-1} c^j x^k 
\label{eq:explicitq}
\end{equation}
\end{subequations}
The binomial coefficients immediately explain the maximum size of the prime divisors of the coefficients of these polynomials.  Furthermore, equation (\ref{eq:explicitpolynomials}) provides a good starting point for investigating any combinatorial meaning of $ P_n $ and  $ Q_n $, which is discussed at the end of this paper.

\section{Fractional Linear Transformations}

Let $ \Cstar $ denote the Riemann sphere $ \mathbb{C} \cup  \{   \infty \} $. A fractional linear transformation is a map $ \beta : \Cstar \onto \Cstar $ of the form
$$
\beta(\tau) = \frac{a\tau+b}{c\tau+d}
$$
where $a,b,c,d \in \mathbb{C}$ and $ ad-bc \neq 0 $.  This function is a conformal map which is analytic everywhere except the pole at $ \tau = -\frac{d}{c} $.  For our purposes, we define a particular fractional linear transformation $ \varphi(\tau) $ as
\begin{equation}
\varphi(\tau) = \frac{\tau - r_1}{\tau - r_2}
\label{eq:phi}
\end{equation}
where $ r_1 $ and $ r_2 $ are the two distinct roots of the quadratic polynomial $ f(x) = ax^2 + bx + c $.  This particular approach exploits the fact that Newton's method for quadratics (with distinct roots) is conjugate to $ z \rightarrow z^2 $, with respect to $ \varphi $, the fractional linear transformation which sends the roots of the quadratic to $ 0 $ and $ \infty $.  Rick Kenyon \cite{kenyon} was the first to point this out to us; expositions can be found in both Cayley \cite{cayley1} and McMullen \cite{mcmullen}.  For distinct roots $ r_1 $ and $ r_2 $, the fractional linear transformation $ \varphi^{-1}(\tau) $ exists and the formula from equation (\ref{eq:newtonrecurrence}) can be expressed as
$$
N(z_n) = z_n - \frac{az_n^2+bz_n+c}{2az_n+b} = \varphi^{-1}(\varphi(z_n)^2)
$$
That is, the following diagram of maps is commutative:
\[
\begin{CD} 
\mathbb{C} @> N >> \mathbb{C} \\
@ V \varphi VV @ VV \varphi V \\
\Cstar @>> (\cdot)^2 > \Cstar \\
\end{CD}
\]
Thus, the formula for the $n$th iterate of Newton's method is
\begin{equation}
\varphi^{-1}(\varphi(x)^{2^n}) = \frac{r_1(x-r_2)^{2^n}-r_2(x-r_1)^{2^n}}{(x-r_2)^{2^n}-(x-r_1)^{2^n}} = \frac{P_n(x)}{Q_n(x)}
\end{equation}
Using this formula, we shall prove the equations in (\ref{eq:explicitpolynomials}).

\section{Proof of Explicit Formula}
\begin{theorem} Given a quadratic polynomial $ f(x) = ax^2 + bx + c $, with $ b^2 - 4ac \neq 0 $, define polynomials $ P_n $ and $ Q_n $ as follows. \begin{subequations}
\begin{equation}
P_n(x) = \frac{a^{2^n-1}}{r_1-r_2} \left( r_1 \left( x-r_2 \right)^{2^n}-r_2 \left( x-r_1 \right )^{2^n}\right) \end{equation}
\begin{equation}
Q_n(x) = \frac{a^{2^n-1}}{r_1-r_2} \left( \left( x-r_2 \right)^{2^n}- \left( x-r_1 \right)^{2^n} \right)
\end{equation}
\label{eq:fltpolynomials}
\end{subequations}
Then the polynomials given in equation (\ref{eq:fltpolynomials}) are equal to those in equation (\ref{eq:explicitpolynomials}) which are reproduced below.
\begin{subequations}
\begin{equation*}
P_n(x) = a^{2^n-1} x^{2^n} - \sum_{k=0}^{2^n-2} ~ \sum_{j=0}^{2^n-k-2} (-1)^j \binom{2^n}{k} \binom{2^n-k-j-2}{j} a^{k+j} b^{2^n-k-2j-2} c^{j+1} x^k
\end{equation*}
\begin{equation*}
Q_n(x) = \sum_{k=0}^{2^n-1}~ \sum_{j=0}^{2^n-k-1} (-1)^j \binom{2^n}{k} \binom{2^n-k-j-1}{j} a^{k+j} b^{2^n-k-2j-1} c^j x^k
\end{equation*}
\end{subequations}
\label{thm:maintheorem}
\end{theorem}
Before we begin the proof of the theorem, we will need the following lemma.
\begin{lemma}
For all $ x, y \in \R $, and all $ n \in \mathbb{N} $, $ n \geq 1 $, the following identity holds
\begin{equation}
x^n - y^n = (x - y) \sum_{i=0}^{n-1} (-1)^i \binom{n-i-1}{i} (x + y)^{n-2i-1} (xy)^i
\label{eq:identity}
\end{equation}
\label{lemma:binomial}
\end{lemma}
\begin{proof}[Proof of Lemma \ref{lemma:binomial}]
We proceed by induction. For $ n = 1, 2 $ the equality is easily verified, so assume the identity is valid for all $ k \leq n $. 
Let $ T(n) $ be the right-hand side of equation (\ref{eq:identity}).  Since
$$
x^{n+1} - y^{n+1} = (x + y) (x^n - y^n ) - xy (x^{n-1} - y^{n-1})
$$
we need only prove that
\begin{equation}
T(n + 1) = (x + y) T(n) - xy T(n-1)
\label{eq:lemmainduction}
\end{equation}
Combining the sums on the right side of equation (\ref{eq:lemmainduction}) by shifting indices, we have
$$
(x - y) \left[ \binom{n-1}{0} (x + y)^{n} + \sum_{i=1}^{n-1} (-1)^{i} \left [ \binom{n-i-1}{i} + \binom{n-i-1}{i-1} \right ] (x + y)^{n-2i} (xy)^i \right ]
$$
Applying Pascal's identity and bringing the leftmost term inside the sum, we obtain 
$$
(x - y) \sum_{i=0}^{n-1} (-1)^{i} \binom{n-i}{i} (x + y)^{n-2i} (xy)^i
$$
Finally, we note that  for $ n \geq 1 $ the binomial coefficient is $ 0 $  when $ i=n $. Therefore, this last expression is equal to $ T(n+1) $, and so we are done.
\end{proof}
\begin{proof}[Proof of Theorem \ref{thm:maintheorem}]
Begin by expanding equation (\ref{eq:fltpolynomials}) via the binomial theorem to obtain
\begin{eqnarray*}
P_n(x)	&=& \frac{a^{2^n-1}}{r_1-r_2} \sum_{k=0}^{2^n} \binom{2^n}{k} x^k \left [ r_1(-r_2)^{2^n-k} - r_2(-r_1)^{2^n-k} \right ] \\
		&=& \frac{a^{2^n-1}}{r_1-r_2} \sum_{k=0}^{2^n} \binom{2^n}{k} x^k \left [ r_1 r_2 \left ((-r_1)^{2^n-k-1} - (-r_2)^{2^n-k-1} \right ) \right ] \\
		&=& \frac{a^{2^n-1}}{r_1-r_2} \left [ (r_1-r_2) x^{2^n} + \sum_{k=0}^{2^n-2} \binom{2^n}{k} x^k \left [ r_1 r_2 \left ((-r_1)^{2^n-k-1} - (-r_2)^{2^n-k-1} \right ) \right ]  \right ] \\
Q_n(x) &=& \frac{a^{2^n-1}}{r_1-r_2} \sum_{k=0}^{2^n-1} \binom{2^n}{k} x^k \left[{(-r_2)^{2^n-k} - (-r_1)^{2^n-k}}\right ] \\
\end{eqnarray*}
We can now apply lemma (\ref{lemma:binomial}) to the expressions of the form $ x^n -y^n $ in both $ P_n $ and $ Q_n $.   Doing so and canceling the factor of $ r_1 - r_2 $ (since $ b^2 - 4ac \neq 0 $) yields
\begin{eqnarray*}
P_n(x) &=& a^{2^n-1} \left [ x^{2^n}  -  \sum_{k=0}^{2^n-2} \binom{2^n}{k} x^k \left [ \sum_{j=0}^{2^n-k-2} (-1)^j \binom{2^n-k-j-2}{j} (-r_1 - r_2)^{2^n-k-2j-2} (r_1 r_2)^{j + 1} \right ] \right ] \\
Q_n(x) &=& a^{2^n-1} \sum_{k=0}^{2^n-1} \binom{2^n}{k} x^k \left [ \sum_{j=0}^{2^n-k-1} (-1)^j \binom{2^n-k-j-1}{j} (-r_1 - r_2)^{2^n-k-2j-1} (r_1 r_2)^j \right ]
\end{eqnarray*}
Replacing $ r_1 $ and $ r_2 $ with their values in terms of the coefficients $ a,b,c $ gives the final form
\begin{eqnarray*}
P_n(x) &=& a^{2^n-1} x^{2^n}  -  \sum_{k=0}^{2^n-2} \binom{2^n}{k} x^k \left [ \sum_{j=0}^{2^n-k-2} (-1)^j \binom{2^n-k-j-2}{j} a^{k+j} b^{2^n-k-2j-2} c^{j+1} \right ] \\
Q_n(x) &=& \sum_{k=0}^{2^n-1} \binom{2^n}{k} x^k \left [ \sum_{j=0}^{2^n-k-1} (-1)^j \binom{2^n-k-j-1}{j} a^{k+j} b^{2^n-k-2j-1} c^j \right ]
\end{eqnarray*}
These equations are the same as equations (\ref{eq:explicitp}) and (\ref{eq:explicitq}), which was to be proved.
\end{proof}


\section{Further Remarks}

Though this paper does not discuss any combinatorial interpretation of the polynomials $ P_n $ and $ Q_n $ we suspect that their may be some fruitful combinatorial equivalence yet to be discovered.  To aid further research in this area, we make a couple of observations about $ P_n $ and $ Q_n $.

Simply iterating Newton's method for the general quadratic gives
\begin{equation}
N \left ( \frac{P_n}{Q_n}  \right ) = \frac{aP_n^2-cQ_n^2}{2aP_nQ_n+bQ_n^2} = \frac{P_{n+1}}{Q_{n+1}}
\label{eq:recurrencepolynomials}
\end{equation}
where $ P_0(x) = x $ and $ Q_0(x) = 1 $ so that the initial term is $ z_0 = x $ as before.  The question is whether the numerator and denominator of equation (\ref{eq:recurrencepolynomials}) are relatively prime so that we may define $ P_n $ and $ Q_n $ recursively in the natural way.  This is, in fact, true as the following lemma proves.
\begin{lemma}
The polynomials $ P_{n+1} $ and $ Q_{n+1} $, defined recursively as
$$ P_{n+1} = aP_n^2-cQ_n^2 $$ and $$ Q_{n+1} = 2aP_nQ_n+bQ_n^2 $$
where $ P_0(x) = x $ and $ Q_0(x) = 1 $, are relatively prime except possibly in the case $ b^2-4ac = 0 $.
\label{eq:relativelyprime}
\end{lemma}
\noindent (For the duration of the statement and proof of Lemma 2, we are suspending 
the definition of $P_n$ and $Q_n$ given earlier, but it will be an
immediate consequence of Lemma 2 that the two definitions agree.)
\begin{proof}
We proceed inductively by assuming that $ P_i $ and $ Q_i $ are relatively prime for all $ i \leq n $.  Assume $ P_{n+1} $ and $
Q_{n+1} $ are not relatively prime to derive a contradiction. Then there exists an irreducible polynomial $ \alpha $ such that $ \alpha~|~P_{n+1} $ and $ \alpha~|~Q_{n+1} $.  If $ \alpha~|~Q_n $, then, since $ \alpha~|~aP_n^2-cQ_n^2 $, it follows that $ \alpha~|~P_n $ which contradicts the induction hypothesis that $ P_n $ and $ Q_n $ are relatively prime.  Hence $ \alpha \not |~Q_n $.  Since $ \alpha~|~2aP_nQ_n+bQ_n^2 $, we know that $ \alpha~|~2aP_n + bQ_n $.  But then $ \alpha~|~P_n + \frac{b}{2a}Q_n $ and also $ \alpha~|~P_n \pm \sqrt{\frac{c}{a}}Q_n $.  Consequently, $ \alpha $ divides their difference, so $ \alpha~|~\left(\frac{b}{2a} \mp \sqrt{\frac{c}{a}} \right)Q_n $ which only occurs when $ b^2 - 4ac = 0$.
\end{proof}

In combinatorics it is sometimes useful to consider two formal variables $ x, y $ which do not commute with each other but instead obey $ yx = qxy $ where $ q $ is another formal variable that commutes with both $ x $ and $ y $.  This approach, due to Sch{\"u}tzenberger \cite{qbin}, is useful in applications such as counting lattice paths.  In our case, the polynomials $ P_n $ and $ Q_n $ can be easily generalized to the non-commuting case.
As a generalization of the usual binomial coefficient, the $ q $-binomial coefficient is defined as 
$$
\qbinomial{n}{k} = \prod_{i=1}^{n-k} \frac{1 - q^{i+k} }{1 - q^i}
$$
and the following $ q $-binomial theorem for non-commuting variables $ x, y $ due to Sch{\"u}tzenberger \cite{qbin} is
$$
(x + y)^n = \sum_{k=0}^{n} \qbinomial{n}{k} x^k y^{n-k} 
$$
Analogously to the above, we can then define non-commuting polynomials $ P_n' $ and $ Q_n' $ recursively as
\begin{subequations}
\begin{equation}
P_{n+1}' = a P_n'^2 - c Q_n'^2
\end{equation}
\begin{equation}
Q_{n+1}' = a P_n' Q_n' + a Q_n' P_n' + b Q_n'^2
\end{equation}
\end{subequations}
where $ P_0' = x $, $ Q_0' = y $, and $ yx = qxy $.  We conjecture explicit formulas for both $ P_n' $ and $ Q_n' $ which happen to be the same as equations (\ref{eq:explicitp}) and (\ref{eq:explicitq}) except for the presence of a $ q $-binomial coefficient:
\begin{subequations}
\begin{equation}
P_n'(x,y) = a^{2^n-1} x^{2^n} - \sum_{k=0}^{2^n-2} ~ \sum_{j=0}^{2^n-k-2} (-1)^j \qbinomial{2^n}{k} \binom{2^n-k-j-2}{j} a^{k+j} b^{2^n-k-2j-2} c^{j+1} x^k y^{2^n - k}
\end{equation}
\begin{equation}
Q_n'(x,y) = \sum_{k=0}^{2^n-1}~ \sum_{j=0}^{2^n-k-1} (-1)^j \qbinomial{2^n}{k} \binom{2^n-k-j-1}{j} a^{k+j} b^{2^n-k-2j-1} c^j x^k y^{2^n - k}
\end{equation}
\end{subequations}

\section{Conclusion}

We have taken the initial observation that Newton's method, when applied to quadratics, produces polynomials with highly composite coefficients and proved an explicit formula for the $ n $th iterate that explains this compositeness as a consequence of the inherent compositeness of binomial coefficients. Furthermore, a recursive definition and a conjectural non-commutative analogue of the polynomials $ P_n $ and $ Q_n $ were noted in hopes of spurring further research into finding a combinatorial interpretation. We believe that a proof of the non-commutative analogue, as well as the larger issue of finding a combinatorial interpretation, are problems which merit further study. It is also worth noting that for higher-degree polynomials, such as cubics, no similar phenomena have been found. In particular, the occurrence of coefficients with large prime factors indicates that no simple product formulas for the coefficients exist, but this does not rule out the existence of more complicated formulas.

We heartily thank Jim Propp, Rick Kenyon, and the rest of the Spatial Systems Laboratory at UW-Madison for their helpful insight and generous support for this research.  In addition, we are indebted to the NSF's Research Experiences for Undergraduates program and the NSA for funding our research, as well as to the computing staff at UW-Madison for providing the computational resources necessary for our investigations.

\bibliographystyle{plain}
\bibliography{newton_bibliography}

\end{document}